\documentclass[12pt]{article}
\usepackage{amsmath}
\usepackage{amssymb}
\usepackage{latexsym}
\def\bbbc{{\mathchoice {\setbox0=\hbox{$\displaystyle\rm
C$}\hbox{\hbox
to0pt{\kern0.4\wd0\vrule height0.9\ht0\hss}\box0}}
{\setbox0=\hbox{$\textstyle\rm C$}\hbox{\hbox
to0pt{\kern0.4\wd0\vrule height0.9\ht0\hss}\box0}}
{\setbox0=\hbox{$\scriptstyle\rm C$}\hbox{\hbox
to0pt{\kern0.4\wd0\vrule height0.9\ht0\hss}\box0}}
{\setbox0=\hbox{$\scriptscriptstyle\rm C$}\hbox{\hbox
to0pt{\kern0.4\wd0\vrule height0.9\ht0\hss}\box0}}}}
\def\bbbe{{\mathchoice {\setbox0=\hbox{\smalletextfont
e}\hbox{\raise
0.1\ht0\hbox to0pt{\kern0.4\wd0\vrule width0.3pt
height0.7\ht0\hss}\box0}}
{\setbox0=\hbox{\smalletextfont e}\hbox{\raise
0.1\ht0\hbox to0pt{\kern0.4\wd0\vrule width0.3pt
height0.7\ht0\hss}\box0}}
{\setbox0=\hbox{\smallescriptfont e}\hbox{\raise
0.1\ht0\hbox to0pt{\kern0.5\wd0\vrule width0.2pt
height0.7\ht0\hss}\box0}}
{\setbox0=\hbox{\smallescriptscriptfont e}\hbox{\raise
0.1\ht0\hbox to0pt{\kern0.4\wd0\vrule width0.2pt
height0.7\ht0\hss}\box0}}}}
\def\bbbr{{\rm I\!R}} 
\def\bbbn{{\rm I\!N}} 

\newcommand{\mn}{\medskip\noindent}
\newcommand{\bn}{\bigskip\noindent}
\newcommand{\sn}{\smallskip\noindent}

\newcommand{\D}{{\cal{D}}}
\newcommand{\A}{{\cal{A}}}

\newcommand{\J}{{\cal{J}}}

\newcommand{\Hh}{{\cal{H}}}

\newcommand{\cD}{{\cal{D}}}

\newcommand{\ii}{{\rm{i}}}

\newcommand{\R}{{\bbbr}}

\newcommand{\llabel}[1]{\label{#1}}
\begin{document}
\begin{center} 
{\LARGE On holomorphic functions on a strip in the complex plane\par}
\end{center}

\centerline{Konrad Schm\"udgen}
\vspace*{1ex}
\centerline{\footnotesize AMS Subject class. (1991): Primary 30D05; Secondary 81S05, 47D40, 17B37.}
\begin{abstract} 
Let $f$ be a holomorphic function on the strip 
$\{z\in\bbbc:-\alpha <{\rm Im}~ z<\alpha\}, 
\alpha>0$, belonging to the 
class $\Hh (\alpha,-\alpha;\varepsilon)$ defined below. It is shown that 
there exist holomorphic functions $w_1$ on $\{z\in\bbbc :0<{\rm Im}~ z<2\alpha\}$ 
and $w_2$ on $\{z\in\bbbc :-2\alpha <{\rm Im}~ z<2\alpha\}$ such that 
$w_1$ and 
$w_2$ have boundary values of modulus one on the real axis and satisfy the 
relation $w_1(z)=f(z{-}\alpha \ii)w_2(z{-}2\alpha \ii)$ and 
$w_2(z{+}2\alpha i)=\bar{f}(z{+}\alpha \ii) w_1(z)$ for $0<{\rm Im}~ z<2$, where 
$\bar{f}(z):={\overline{f(\bar{z})}}$. This leads to 
a "polar decomposition" $f(z)=u_f(z{+}\alpha \ii)g_f(z)$ of the 
function $f(z)$, where $u_f(z{+} \alpha\ii)$ and $g_f(z)$ are holomorphic 
functions for $-\alpha<{\rm Im}~ z<\alpha$ such that $|u_f(x)|=1$ and 
$g_f(x)\ge 0$ a.e. on the real axis. As a byproduct, an operator 
representation of a $q$-deformed Heisenberg algebra is developed.
\end{abstract}
\renewcommand{\baselinestretch}{1.0}
\section{Introduction and Main Results}

Let $\epsilon$ be a positive number and $\alpha$
and $\beta$ be real numbers such that $\alpha > \beta$. Let 
${\cal H}(\alpha,\beta;\epsilon)$ denote the set of all holomorphic functions 
$h(z)$ on the strip ${\cal I}(\alpha,\beta) 
:=\{z\in\bbbc:\alpha >  {\rm Im}~ z >\beta\}$ such that
$$
\sup_{\alpha > y > \beta}\,\,\int\limits_{-\infty}^\infty\,
\vert h\,(x{+}y\ii)\vert^2\,e^{-2\gamma x^2}\,dx\;<\;\infty
$$
for all numbers $\gamma >\epsilon$. As stated in Lemma 2 below, each
function $h\in{\cal H}(\alpha,\beta;\epsilon)$ admits boundary values 
$h(x{+}\beta \ii)$
and $h(x{+}\alpha \ii),x\in\bbbr$, which satisfy
\begin{eqnarray}\label{1}
\lim_{y\downarrow \beta}\,\,\int\,\vert h\,(x{+}y\ii) {-} 
h\,(x{+}\beta \ii)\vert^2\,
e^{-2\gamma x^2}\,dx~ =~ 0,\\
\label{2} 
\lim_{y\uparrow \alpha}\,\,\int\,\vert h\,(x{+}y\ii) {-} h\,(x{+}\alpha\ii)
\vert^2\,e^{-2\gamma x^2}\,dx~=~0 
\end{eqnarray}
for all $\gamma >\epsilon$. Throughout this paper, 
${\rm i}$ denotes the complex unit. 
By some slight abuse of notation, we  
denote functions in ${\cal H}(\alpha,\beta;\epsilon)$ and their boundary 
values by the same symbol.

\bigskip
Our main results are contained in the following

\bigskip\noindent
{\bf Theorem 1.} {\it Let $\varepsilon$ and $\alpha$ be positive 
real numbers and let $f\ne 0$ a function of the class 
$\Hh (\alpha, -\alpha;\varepsilon)$ such that
\begin{equation}\label{b1}
\inf \{ |f (x{-}\alpha\ii)|; x\in\bbbr\} > 0~.
\end{equation}
Then there exist functions $w_1\in\Hh (2\alpha,0;\varepsilon)$ and  
$w_2\in \Hh(2\alpha,-2\alpha;\varepsilon)$ such that $|w_1 (x)|=|w_2 (x)|=1$ a.e. on $\bbbr$ and
\begin{align}\label{b2}
w_1(x)&=f(x{-}\alpha \ii) w_2(x{-}2\alpha \ii),\\
\label{b3}
w_2(x)&=\bar{f}(x{-}\alpha \ii)w_1(x{-}2\alpha \ii)
\end{align}
a.e. on $\bbbr$. If $\tilde{w_1}$ and $\tilde{w_2}$ are two other 
functions with these properties, then there is a constant $c$ of modulus 
one such that $\tilde{w_1} (x)= c w_1(x)$ and 
$\tilde{w_2} (x) = c w_2 (x)$ a.e. on $\bbbr$.}

\bigskip
We briefly discuss some consequences of the preceding result. Since both 
sides of the equality (\ref{b2}) are boundary values of holomorphic 
functions on the strip $\J(2\alpha,0)$, we conclude that
\begin{equation}\label{b4}
w_1(z{+}2\alpha \ii)=f(z{+}\alpha \ii) w_2(z),~ z\in \J(0,{-}2\alpha).
\end{equation}
Since $\bar{f}(z{-}\alpha \ii)w_1(z{-}2\alpha \ii)$ is holomorphic on 
$\J(2\alpha, 0)$, it follows from (\ref{b3})
that the function $w_2(x)$ on $\bbbr$ is boundary value of a holomorphic 
function  $w_2(z)$ on $\J(2\alpha, 0)$ and 
\begin{equation}\label{b5}
w_2 (z{+}2\alpha \ii)= \bar{f}(z{+}\alpha \ii) w_1 (z),~ z\in\J(2\alpha,0).
\end{equation}

If the function $f\in\Hh(\alpha,-\alpha;\varepsilon)$ from Theorem 1 is 
holomorphic on the whole upper half-plane, then it follows from (\ref{b4}) 
and (\ref{b5}) that $w_1(z)$ and $w_2(z)$ are  holomorphic on the upper 
half-plane and that the relations (\ref{b4}) and (\ref{b5}) hold for all 
$z\in\bbbc, {\rm Im}~ z>-2\alpha$. 
Moreover, both relations then imply that
\begin{align*}
w_1(z{+}4\alpha \ii)&=f(z{+}3\alpha \ii)\bar{f}(z{+}\alpha \ii) w_1(z),\\
w_2(z{+}4\alpha \ii)&=\bar{f}(z{+}3\alpha \ii)f(z{+}\alpha \ii) w_2(z)
\end{align*}
for ${\rm Im}~z>-2\alpha$.

Theorem 1 is a result on holomorphic functions, but its proof is based on 
operator-theoretic tools. This technique might be also of interest in itself.
 Let $P$ be the operator $-i{\frac{d}{dx}}$ on $L^2(\bbbr)$. 
 For $f\in\Hh(\alpha,-\alpha;\varepsilon)$, let $L_f$ and $R_f$ 
 be the closures of the linear operators $\tilde{L}_f$ and 
 $\tilde{R}_f$ on the Hilbert space $L^2(\bbbr)$ defined by the formulas
$$
\tilde{L}_f=f(x{-}\alpha \ii)e^{2\alpha P}~{\rm and}~ 
\tilde{R}_f=e^{2\alpha P} \bar{f}(x{+}\alpha\ii).
$$
The operators $L_f$ and $R_f$ are crucial in the proof of 
Theorem 1. The polar decomposition of the operator $L_f$ 
is described by the following theorem. Some more properties of these 
operators can be found in Sections 2 and 3 and in our previous 
papers [5] and [6].

\bigskip\noindent
{\bf Theorem 2.} {\it Retain the assumptions and notations from Theorem 1 
and define holomorphic functions $u_f$ and $g_f$ 
on the strips 
$\J(2\alpha,0)$ and $\J(\alpha,-\alpha)$, respectively, by}
\begin{equation}\label{b6}
u_f(z)=w_1(z)\overline{w_2}(z), ~g_f(z)=
w_2(z{+}\alpha \ii)\overline{w_2}(z{-}\alpha \ii).
\end{equation}
{\it Then the polar decomposition of the 
closed linear operator $L_f$ is given by $L_f=u_f L_{g_f}$ and 
we have $L_f^\ast=R_f$.}

\bigskip
From the relation $L_f=u_fL_{g_f}$ we immediately obtain that 
$f(x{-}\alpha \ii)=u_f(x) g_f(x{-}\alpha \ii)$ a.e. on $\bbbr$. This in turn 
implies that
\begin{equation}\label{b7}
f(z)=u_f(z{+}\alpha \ii) g_f(z),~ z\in\J(\alpha,-\alpha),
\end{equation}
where $u_f(z{+} \alpha \ii)$ and $g_f(z)$ are holomorphic functions on the 
strip $\J(\alpha,-\alpha)$ such that
\begin{align*}
&|u_f(x)|=1 \qquad\text{a.e. on}~\bbbr,\\
&g_f(x) = |w_2 (x{+}\alpha \ii)|^2\ge 0~{\rm on}~\bbbr.
\end{align*}
Because of the two latter properties and the fact that it originates from 
the polar decomposition of the operator $L_f$, we refer to the decomposition 
(\ref{b7}) as the {\it polar decomposition} of the holomorphic function 
$f\in\Hh(\alpha,-\alpha;\varepsilon)$ with respect to the strip 
$\J(\alpha,-\alpha)$.

The poroofs of Theorems 1 and 2 will be given in Section 3. Some 
necessary operator-theoretic tools will be collected in Section 2. 
The example $f(z)=z$ is treated in Section 4. As an application we 
construct in Section 5 an interesting operator representation of a 
$q$-deformed Heisenberg algebra.  



\section*{2\ Technical Preliminaries}

\bigskip
Throughout, $P$ denotes the self-adjoint operator $-i{d\over dx}$ 
on the Hilbert space $L^2(\bbbr)$ and $\alpha$ is a positive real number. 
The following lemma describes the action and the domain of the operator 
$e^{\alpha P}$. Its proof given in [6] is 
essentially based on the classical Paley-Wiener theorem [2]. The domain of 
an operator $T$ is denoted by ${\cal D}(T)$.

\bigskip\noindent
{\bf Lemma 1.} {\it Suppose that $\alpha > 0$. Let $g(z)$ be a
holomorphic function on the strip ${\cal I}(0,-\alpha)$ such that
$$
\sup_{0<y<\alpha}\;\int\limits_{-\infty}^\infty\;\vert g(x{-}iy)\vert^2\,
dx\;<\;\infty\;.\eqno(6)
$$
Then there exist functions $g(x)\in L^2(\bbbr)$ and $g_{-\alpha}(x)\in
L^2(\bbbr)$ such that $\lim_{y\uparrow 0}\,g_y=g$ and 
$\lim_{y\downarrow -\alpha}\,g_y=g_{-\alpha}$ in $L^2(\bbbr)$, where 
$g_y(x):=g(x{+}iy)$ for $x\in\bbbr$ and $y\in(0,-\alpha)$. Setting 
$g(x{-}i\alpha):=g_{-\alpha}(x), x\in\bbbr$, we have
$$
\lim_{n\to\infty}\,g\left(x{-}n^{-2}\ii \right)=g(x)~{\rm and}~
\lim_{n\to\infty}g\left(x{-}(\alpha{+}n^{-2})\ii\right)=
g(x{-}\alpha \ii)~{\rm a.e.~on}~ \bbbr\;.\eqno(7)
$$
The function $g$ belongs to the domain ${\cal D}(e^{\alpha P})$ and 
$(e^{\alpha P} g)(x)=g(x{-}\alpha\ii)$. Conversely, each function  
$g\in{\cal D}(e^{\alpha P})$ arises in this way.}

\bigskip\noindent
{\bf Proof.} [6], Lemma 1. \hfill $\Box$

\bigskip
For $\delta \geq 0$, we define a dense linear subspace 
$\cD_\delta$ of the Hilbert space $L^2(\bbbr)$ by
\begin{equation}\label{c9} 
{{\cal D}_\delta =
{\rm Lin}\{e^{-\gamma x^2{+}\beta x}{:}} 
~~\gamma>\delta, \beta\in\bbbc\}.
\end{equation}

\bigskip\noindent
{\bf Lemma 2.} (i) {\it Suppose that $h\in{\cal H}(\alpha,\beta ;\epsilon)$.
 Then there exist
measurable functions $h(x{+}\beta \ii)$ and $h(x{+}\alpha \ii)$ on 
$\bbbr$ both
contained in the domain ${\cal D}(e^{-\gamma x^2})$ for each 
$\gamma >\epsilon$ such that (\ref{1}) and (\ref{2}) 
are satisfied. If $\varphi \in L^2(\R)$ such that 
$h\varphi \in {\cal D}(e^{\alpha P})$ (in particular, 
if $\varphi \in \cD_\varepsilon$), we have}
\begin{eqnarray}\label{c10}
\left(e^{\alpha P}h\varphi\right)(x)=
h(x{-}\alpha\ii)\varphi(x{-}\alpha\ii)=h(x{-}\alpha\ii)
\left(e^{\alpha P}\varphi\right)
(x).
\end{eqnarray} 
(ii) {\it For any $\alpha\in\bbbr$ and $\beta \geq 0$, ${\cal D}_\beta$ is a
core for $e^{\alpha P}$.}

\bigskip\noindent
{\bf Proof.} [6], Lemma 2 and Lemma 3(ii). \hfill$\Box$

\bigskip\noindent
{\bf Lemma 3.}{\it Let $\alpha > 0$. Suppose that $f\ne 0$ is a 
function of the class $\Hh (\alpha, -\alpha;\varepsilon)$ 
satisfying condition (\ref{b1}). Then we have:}\\
(i) $L_f = {\tilde L}_f \equiv  f(x{-}\alpha \ii) e^{2\alpha P}.$\\
(ii) $\ker\,L_f =\{0\}$.\\
(iii) $L_f^\ast = R_f.$

\mn
{\bf Proof.} (i): We have to show that the operator 
${\tilde L}_f = f(x{-}\alpha \ii) e^{2\alpha P}$ is closed. 
For let $\{\eta_n \}$ be a sequence of vectors $\eta_n \in 
\cD({\tilde L}_f)$ such that $\eta_n \to \eta$ and 
$f(x{-}\alpha \ii) e^{2\alpha P}\eta_n \to \zeta$ in the Hilbert 
space $L^2(\R)$. Since $\{f(x{-}\alpha \ii) e^{2\alpha P}\eta_n\}$ is a 
Cauchy sequence, condition (\ref{b1}) implies that 
$\{e^{2\alpha P}\eta_n\}$ is also a Cauchy sequence in $L^2(\R)$, so that 
$e^{2\alpha P}\eta_n \to \xi$ for some $\xi \in L^2(\R)$. Because the operator $e^{2\alpha P}$ is closed, we 
conclude that $\xi \in \cD(e^{2\alpha P})$ and $e^{2\alpha P}\eta = \xi$. 
Since the multiplication operator by the function $f(x{-}\alpha \ii)$ is 
also closed, $\eta \in \cD(f(x{-}\alpha \ii) e^{2 \alpha P})$ 
and ${\tilde L}_f \eta \equiv f(x{-}\alpha \ii) e^{2 \alpha P} \eta = 
f(x{-}\alpha \ii)\xi = \zeta $. \\
(ii) follows immediately from (i).\\
(iii): Using the relation ${\overline{f(x{-}\alpha \ii)}} = 
{\overline f}(x{+}\alpha \ii)$ one easily verifies that 
$R_f \subseteq L_f^\ast$. We now prove the opposite inclusion. Since the 
function ${\overline f}(x{+}\alpha \ii)^{-1}$ is bounded because of 
assumption (\ref{b1}), we have $e^{2 \alpha P} \subseteq R_f {\overline f}(x{+}\alpha \ii)^{-1}$. Using this fact we conlcude that 
$$
f(x{-}\alpha \ii)^{-1} R_f^\ast \subseteq (R_f {\overline f}(x{+}\alpha \ii)^{-1})^\ast \subseteq (e^{2 \alpha P})^\ast = e^{2 \alpha P}
$$
which in turn implies that 
$$ 
R_f^\ast \subseteq f(x{-}\alpha \ii)e^{2 \alpha P} = L_f.
$$
Applying the adjoint to the latter we get $L_f^\ast \subseteq 
(R_f)^{\ast\ast} = R_f$, because the operator $R_f$ is closed. 
Both inclusion together yield the desired equality 
$L_f^\ast = R_f.$\hfill$\Box$

\bn
\section*{3\ Proof of the Theorems}
Consider the self-adjoint operator $\A$ and the one parameter unitary 
group $U(t), t\in\R$, on the Hilbert space $\Hh=L^2(\R)\oplus L^2(\R)$
 given by the operator matrices
\begin{eqnarray}\label{af}
A_f=\left(\begin{matrix}0 &L_f\\ L_f^\ast &0\end{matrix}\right),
\qquad U(t)=\left(\begin{matrix} e^{\ii t Q} & 0\\ 0 &e^{\ii t Q}
\end{matrix}\right).
\end{eqnarray}
Since $e^{-\ii t Q}~e^{2\alpha P} ~e^{\ii t Q}= e^{2\alpha t+2\alpha P}$,
 we have $e^{-\ii t Q} L_f e^{\ii t Q}=e^{2\alpha t} L_f$ and hence 
$e^{-\ii tQ} L_f^\ast e^{\ii tQ}= e^{2\alpha t} L_f^\ast$ for $t\in\bbbr$. 
These relations immediately imply that
\begin{equation}\llabel{a1}
U(-t)A_f U(t)=e^{2\alpha t} A_f, t\in\R~.
\end{equation}
Let us write $A_f=A^+_f\oplus ({-}A^-_f)$, where $A^\pm_f$ and 
$A^-_f$ are the positive and the negative parts of the self-adjoint 
operator $A_f$. The corresponding reducing subspaces of the Hilbert 
space $\Hh$ are denoted by $\Hh_+$ and $\Hh_-$. We first verify that 
$A^+_f\ne 0$ and $A^-_f\ne 0$. Let $L_f= u_f|L_f|$ be the polar 
decomposition of the closed operator $L_f$. Since $\ker~ L_f=\{0\}$ 
and $\ker~ L_f^\ast=\{0\}$ as noted above, $u_f
$ is unitary. We have $L_f^\ast=|L_f|u_f^\ast$ and hence
\begin{equation}\label{a2}\begin{split}
\langle A_f (\eta_1\oplus \eta_2), \eta_1\oplus \eta_2\rangle 
&=\langle L_f \eta_2, \eta_1\rangle + \langle L_f^\ast \eta_1,\eta_2\rangle\\
&=\langle |L_f|\eta_2, u_f^\ast\eta_1\rangle + 
\langle |L_f|u_f^\ast\eta_1,\eta_2\rangle
\end{split}\end{equation}
for $\eta_1\in\cD (L_f^\ast)$ and $\eta_2\in\cD(L_f)$. 
Choosing $\eta_2=u_f^\ast\eta_1, \eta_1\ne 0$, the expression (\ref{a2}) 
becomes positive. For $\eta_2=-u_f^\ast\eta_1,\eta_1\ne 0$, it is negative.
 This implies that $A^+_f\ne 0$ and $A^-_f\ne 0$.

From relation (\ref{a1}) it follows that $U(-t) A^\pm_f U(t)=e^{2\alpha t} A^\pm_f$ 
and  $U(t)\Hh_\pm\subseteq \Hh_\pm$ for real $t$. Recall that $A^+_f$ and 
$A^-_f$ are positive self-adjoint operators with trivial kernels. Therefore, 
we conclude that $U(-t) (A^\pm_f)^{s\ii} U(t)
=(e^{2\alpha t} A^\pm_f)^{s\ii}=
e^{2\alpha st \ii} (A^\pm_f)^{s\ii}$ for $s,t\in\bbbr$. That is, the unitary 
groups $U_\pm(t):=U(t)\lceil \Hh_\pm$ and $V
_\pm(s):=(A^\pm_f)^{is}$ on the Hilbert space $\Hh_\pm$ satisfy the Weyl 
relation $V_\pm(s) U_\pm(t)=e^{2\alpha st i}U_\pm(t)V_\pm(s),s,t\in\R$. By 
construction, the unitary group $U(t)=U_+(t)\oplus U_-
(t)$ has uniform spectral 
multiplicity two. Consequently, since $A^\pm_f\ne 0$ and hence 
$\Hh_\pm\ne\{0\}$, both unitary groups $U_+(t)$ and $U_-(t)$ have 
spectral multiplicity one. Therefore, it follows from the Stone--von 
Neumann uniqueness theorem (see, for instance, [3]) that each pair 
$\{U_\varepsilon(t), A^\varepsilon_f\}$ on 
$\Hh_\varepsilon, \varepsilon =\pm 1$, is unitarily equivalent to the 
pair $\{ e^{\ii t Q}, e^{2\alpha P}\}$ acting on the Hilbert space $L^2(\R)$. 
Hence the pair $\{ U(t)=U_+(t)\oplus U_-(t), A_f=A^+_f\oplus (-A^-
_f)\}$ is unitarily equivalent to the pair $\{e^{\ii t Q}\oplus e^{\ii t Q}, 
e^{2\alpha P}\oplus (-e^{2\alpha P})\}$ on $L^2(\R)\oplus L^2(\R)$. For 
the subsequent considerations it is convenient to transform the latter pair 
by means of the unitary symmetry
$$
{{1}\over{\sqrt{2}}} \left( \begin{matrix}I &~~I\cr I &-
I\cr\end{matrix}\right)~.
$$
Putting the preceding together, it follows that there exists a unitary 
2$\times$2-operator matrix 
$$
W=\left(\begin{matrix} w_1 & w_3\\ w_4 &w_2\end{matrix}\right)
$$
of the Hilbert space $\Hh=L^2(\R)\oplus L^2(\R)$ such that
\begin{eqnarray}\label{a3}
&W^\ast U(t)W=\left( \begin{matrix} e^{\ii t Q} &0\cr 0 &e^{\ii t Q}\cr 
\end{matrix} \right),~
t\in \R~,\\
\label{a4}
&W^\ast A_f W=\left(\begin{matrix} 0 &e^{2\alpha P}\cr e^{2\alpha P} 
&0\cr\end{matrix}\right) =: B.\qquad\hskip0.2cm
\end{eqnarray}
Relation (\ref{a3}) means that $w_{j} e^{\ii t Q}=e^{\ii t Q} w_{j}, t\in\R$, for 
$j=1,2,3,4$. This implies that the entries $w_{j}$ of the operator 
matrix $W$ are multiplication operators by essentially bounded 
measurable functions $w_{j}(x),~x\in\R$. Equation (\ref{a4}) is equivalent 
to the two relations $WB \subseteq A_fW$ and $W^\ast A_f \subseteq BW^\ast$.
Applying the relation $WB \subseteq A_fW$ to vectors $(0,\eta)$, 
and $(\eta,0)$, where $\eta \in \D(e^{2 \alpha P})$, in the domain of $B$ 
we obtain that
\begin{eqnarray}\label{a5}
&w_{1} e^{2\alpha P}\subseteq L_f w_{2},
~ w_{4}e^{2\alpha P} \subseteq L_f^\ast w_{3},\\
\label{a6}
&w_{3}  e^{2\alpha P}\subseteq L_f w_{4},
~ w_{2} e^{2\alpha P}\subseteq L_f^\ast w_{1},
\end{eqnarray}
respectively. Similarly, the relation $W^\ast A_f \subseteq BW^\ast$ 
applied to vectors  $(0,\eta)$, $\eta \in \D(L_f)$, 
and $(\eta,0)$, $\eta \in \D(L_f^\ast)$, yields 
\begin{eqnarray}\label{a7}
&\overline{w_{1}} L_f \subseteq   e^{2\alpha P} \overline{w_{2}},
~ \overline{w_{3}} L_f  \subseteq e^{2\alpha P} \overline{w_{4}},\\
\label{a8}
&\overline{w_{4}}  L_f^\ast \subseteq e^{2\alpha P} \overline{w_{3}},
~ \overline{w_{2}} L_f^\ast \subseteq e^{2\alpha P} \overline{w_{1}},
\end{eqnarray}
respectively. In fact, (\ref{a5}) and (\ref{a6}) are equivalent to the inclusion 
$WB \subseteq A_fW$, while (\ref{a7}) and (\ref{a8}) are equivalent to 
$W^\ast A_f \subseteq BW^\ast$. 

From (\ref{a5}) and (\ref{a6}) it follows at once that 
\begin{eqnarray}\label{a9}
&w_{1} e^{4\alpha P}\subseteq L_f L_f^\ast w_{1},~ w_{2} e^{4\alpha P} 
\subseteq L_f^\ast L_f w_{2},\\
\label{a10}
&w_{3} e^{4\alpha P} \subseteq L_f L_f^\ast w_{3},~ w_{4} e^{4\alpha P} 
\subseteq L_f^\ast L_f w_{4}~.
\end{eqnarray}
Let us fix one of the relations $w B_1\subseteq B_2 w$ of 
(\ref{a9}) or (\ref{a10}), where $w=w_{j}, B_1=e^{4\alpha P}$ and $B_2$ is 
one of the self-adjoint operator $L_fL_f^\ast$ or $L_f^\ast L_f$, 
respectively. Since $B_1$ and $B_2$ are self-adjoint and $w$ is bounded, 
we obtain $B_1 w^\ast=(wB_1)^\ast\supseteq (B_2w)^\ast\supseteq w^\ast B_2$ 
and hence $w^\ast w B_1\subseteq w^\ast B_2 w\subseteq B_1 w^\ast w$. 
That is, we have $|w_{j}(x)|^2 e^{4\alpha P}\subseteq e^{4\alpha P} 
|w_{j} (x)|^2$. From the latter we conclude that the bounded operator 
$|w_{j}(x)|^2$ commutes with all functions of the unbounded self-adjoint 
operator $e^{4\alpha P}$, so $|w_{j}(x)|^2$ commutes in particular with 
the unitary group $e^{\ii s P}, s\in\R$, on $L^2(\R)$. Therefore, the 
function $|w_{j} (x)|$ is almost everywhere constant on $\R$, 
say $|w_{j}(x)|=c_{j}$ a.e. on $\R$. 

Since the 2$\times$2-matrix $W$ is a unitary operator on $\Hh$, we 
conclude that $c_1=c_2, c_3=c_4$, and $c_1^2+c_3^2=c_2^2+c_4^2=1$. 
Without loss of generality let us suppose that  $c_1\ne 0$. (If $c_1=0$, 
then $c_3\ne 0$ and we replace the function $w_1,w_2$  by $w_3,w_4$ 
in what follows.) Then, upon replacing $w_1$ by $w_1 c^{-1}_1$ and 
$w_2$ by $w_2 c^{-1}_1$ we can assume that the functions 
$w_1 (x)$ and $w_2(x)$ satisfying (\ref{a5}) and (\ref{a6}) 
are of modulus one on the real line.

From the first relations of (\ref{a5}) and (\ref{a7}) and the second 
relations of (\ref{a6}) and (\ref{a8}) we now easily 
obtain
\begin{eqnarray}\label{a11}
L_f = w_1 e^{2\alpha P} \overline{w_2} ~~{\rm and}~~
L_f^\ast = w_2 e^{2\alpha P} \overline{w_1}.
\end{eqnarray}
Thus, we have $L_f^\ast L_f = w_2 e^{4\alpha P}\overline{w_2}$ 
and $L_f L_f^\ast = w_1 e^{4\alpha P}\overline{w_1}$. For the operators 
$|L_f|=(L_f^\ast L_f)^{1/2}$ and $|L_f^\ast|=(L_f L_f^\ast)^{1/2}$ we 
therefore obtain 
\begin{eqnarray}\label{a12}
|L_f|=w_2 e^{2\alpha P} \overline{w_2} ~~{\rm and}~~
|L_f^\ast|=w_1 e^{2\alpha P} \overline{w_1}.
\end{eqnarray}

\bigskip\noindent
{\bf Remarks.} 1. Recall that $|w_1(x)|=|w_2 (x)|=1$ a.e. on $\bbbr$. 
Using this fact and repeating the above reasoning it follows that 
the two relations (\ref{a11}) are equivalent to
\begin{eqnarray}\label{a13}
\left(\begin{matrix}{\overline w_1} &0\\ 0 &{\overline w_2}\end{matrix}\right) 
\left(\begin{matrix}0&L_f\\L_f^\ast&0\end{matrix}\right)
\left(\begin{matrix}w_1 &0\\ 0&w_2\end{matrix}\right)\quad=\quad 
\left(\begin{matrix}0&e^{2\alpha P}\\ e^{2\alpha P}&0\end{matrix}\right)~.
\end{eqnarray}
That is, the unitary matrix $W$ satisfying (\ref{a3}) and (\ref{a4}) 
can be chosen to be diagonal with functions $w_1$ and $w_2$ as 
diagonal entries.

2. Note that in equations (\ref{a11}) and (\ref{a12}) we have strict 
equality of the corresponding unbounded operators on both sides. 

\mn
Let us begin with the proof of Theorem 1. By Lemma 3(i), 
the operator $\tilde{L}_f=f(x{-}\alpha \ii)e^{2\alpha P}$ is closed 
because of 
assumption (\ref{b1}), so it coincides with $L_f$. From the first 
equality of (\ref{a11}) we obtain
\begin{eqnarray}\label{a14}
f(x{-}\alpha \ii)e^{2\alpha P}=w_1 e^{2\alpha P}\overline{w_2} \qquad \text
{and}\quad f(x{-}\alpha \ii) e^{2\alpha P} w_2=w_1 e^{2\alpha P} .
\end{eqnarray}
Set $\eta_\gamma (x):=e^{-\gamma x^2}$ for $\gamma > \varepsilon$. 
Since $f\in\Hh(\alpha,-\alpha;\varepsilon)$ by assumption, $\eta_\gamma \in \cD (f(x{-}\alpha \ii)e^{2\alpha P})$. Therefore, by the first 
relation of (\ref{a14}), $\overline{w_2} \eta_\gamma\in\cD (e^{2\alpha P})$. 
Since $\eta_\gamma$ is also in the domain of $w_1e^{2\alpha P}$, 
the second equality of (\ref{a14}) yields 
$w_2\eta_\gamma\in\cD(e^{2\alpha P})$. By the characterization of the 
domain $\cD(e^{2\alpha P})$ given in Lemma 1, the facts that 
$\overline{w_2}\eta_\gamma$ and $w_2\eta_\gamma$ for arbitrary 
$\gamma > \varepsilon$ are in $\cD(e^{2\alpha P})$ 
imply that $\overline{w_2}\in\Hh (0,-2\alpha;\varepsilon)$ and 
$w_2\in\Hh(0,-2\alpha;\varepsilon)$. Obviously, the fact that 
$\overline{w_2}\in\Hh(0,-2\alpha;\varepsilon)$ leads to 
$w_2\in\Hh(2\alpha,0;\varepsilon)$. Since $w_2$ is holomorphic on the union 
$\J(2\alpha,0)\cup \J(0,2\alpha)$ and has boundary values of modulus one 
on the real axis, it follows from Morera$^\prime$s theorem that $w_2$ is 
holomorphic on the strip $\J(2\alpha,-2\alpha)$. Having this it is 
clear that $w_2\in\Hh(2\alpha,-2\alpha;\varepsilon)$. Applying 
the second equality of (\ref{a14}) to the vector $\eta_\gamma$ and using  
(\ref{c10}), we get 
$$
f(x{-}\alpha \ii)w_2(x{-}2\alpha\ii)\eta_\gamma (x{-}2\alpha\ii)=
w_1(x)\eta_\gamma(x{-}2\alpha\ii).
$$
Since $\eta_\gamma (x{-}2\alpha\ii)\ne 0$ on $\bbbr$, this gives 
equation ({\ref{b2}).

Next we prove formula (\ref{b3}) and the fact that 
$w_1\in\Hh(2\alpha,0;\varepsilon)$. Combining Lemma 3(iii) and the 
second equality of (\ref{a11}) we get
\begin{equation}\label{a15}
e^{2\alpha P}\overline{f}(x{+}\alpha\ii)\subseteq R_f=L_f^\ast =
w_2 e^{2\alpha P}\overline{w_1}. 
\end{equation}
Note that $\overline{f}\in\Hh(\alpha,-\alpha;\varepsilon)$ because of the 
assumption $f\in\Hh(\alpha,-\alpha;\varepsilon)$. Therefore, 
$\eta_\gamma$ is in the domain of the operator 
$e^{2\alpha P}\overline{f}(x{+}\alpha\ii)$ and hence
$\overline{w_1}\eta_\gamma\in\cD(e^{2\alpha P})$ by (\ref{a15}) for any 
$\gamma > \varepsilon$. From the latter fact we conclude that 
$\overline{w_1}\in\Hh(0,-2\alpha;\varepsilon)$ and so 
$w_1\in\Hh(2\alpha,0;\varepsilon)$. Applying (\ref{a15}) to the vector 
$\eta_\gamma$, we obtain 
$$
\overline{f}(x{-}\alpha\ii)\eta_\gamma(x{-}2\alpha\ii)=
w_2(x)\overline{w_1}(x{-}2\alpha\ii)\eta_\gamma (x{-}2\alpha\ii).
$$
Since $|w_1(x)|=1$ a.e. on $\bbbr$ and hence 
$w_1(z)=\overline{w_1} (z)^{-1}$ for $z\in \J(0,-2\alpha)$, the 
preceding equation implies (\ref{b3}).

Finally, we prove the uniqueness assertion of Theorem 1. Suppose that 
$ \tilde{w}_1$ and $\tilde{w}_2$ are two other functions having the 
properties of $w_1$ and $w_2$, respectively, stated Theorem 1. 
The crucial step of this part of the proof is to show that then
\begin{equation}\label{a16}
\tilde{w}_1 e^{2\alpha P}\subseteq L_f \tilde{w}_2 
\quad\text{and}\quad \tilde{w}_2 e^{2\alpha P}\subseteq  L_f^\ast \tilde{w}_1.
\end{equation}
Let $\eta \in\cD_\varepsilon$. 
Since $\tilde{w}_2\in\Hh(0,-2\alpha;\varepsilon)$ 
by assumption, we have $\tilde{w}_2\eta\in\cD(e^{2\alpha P})$ 
by Lemma 1 and $(e^{2\alpha P}w_2\eta)(x)=
\tilde{w}_2(x{-}2\alpha\ii)\eta (x{-}2\alpha\ii)$ by (\ref{c10}). 
From the relation ${\tilde w}_1(x)=f(x{-}\alpha \ii) 
{\tilde w}_2(x{-}2\alpha \ii)$ by (\ref{b2}) and the preceding equality we 
conclude that $\tilde{w}_2\eta\in\cD(L_f)$ 
and $\tilde{w}_1 e^{2\alpha P}\eta = L_f\tilde{w}_2\eta$. 
By Lemma 2(ii), $\cD_\varepsilon$ is a core for the operator 
$e^{2\alpha P}$. Therefore, since $T_f$ is a closed operator and 
$\tilde{w}_1$ and $\tilde{w}_2$ are unitaries on the Hilbert space 
$L^2(\bbbr)$, the relation $\tilde{w}_1 e^{2\alpha P}\eta = 
LT_f\tilde{w}_2\eta,~\eta\in\cD_\varepsilon$, implies that
$\tilde{w}_1 e^{2\alpha P}\subseteq T_f\tilde{w}_2.$

Because $\tilde{w}_2(z{+}2\alpha\ii),\overline{f}(z{+}\alpha\ii)$ and 
$\tilde{w}_1(z)$ are holomorphic on the strip 
$\J(0,{-}2\alpha)$ by assumption, it follows from (\ref{b3}) that 
\begin{eqnarray}\label{a17}
\tilde{w}_2(x{+}2\alpha\ii)=\overline{f}(x{+}\alpha\ii)\tilde{w}_1(x)
~~{\rm a.e. on}~ \bbbr.
\end{eqnarray}
Since $\tilde{w}_2\in\Hh(2\alpha,0;\varepsilon)$, 
$\tilde{w}_2\eta_\gamma\in\cD(e^{\alpha P})$ for 
$\eta_\gamma\in\cD_\varepsilon$. Therefore, using formulas (\ref{c10}) and 
(\ref{a17}) 
and the fact $L_f^\ast = R_f$ we obtain
\begin{align*}
\tilde{w}_2 (x) e^{2\alpha P}\eta 
= e ^{2\alpha P} \tilde{w}_2 (x{+}2 \alpha \ii) \eta
= e^{2\alpha P} \overline{f} (x {+}\alpha \ii) \tilde{w}_1 (x) \eta
= R_f \tilde{w}_1 \eta = L^\ast_f \tilde{w}_1 \eta.
\end{align*}
Arguing as in the preceding paragraph, the latter implies that 
$\tilde{w}_2e^{2\alpha P}\subseteq L_f^\ast\tilde{w}_1$ which 
proves the second inclusion of (\ref{a16}).

Let $\tilde{W}$ denote the 2$\times$2 diagonal matrix which has 
$\tilde{w}_1$ and $\tilde{w}_2$ in the main diagonals.  Note that 
(\ref{a16}) is nothing but relations (\ref{a5}) and (\ref{a6}) with 
$w_3=w_4=0$ and $w_j$ replaced by $\tilde{w}_j,j=1,2$. Therefore, 
as noted after formula (\ref{a8}), relations (\ref{a16}) are equivalent 
to the inclusion ${\tilde W}B \subseteq A_f{\tilde W}$ and so 
to ${\tilde W}B{\tilde W}^\ast \subseteq A_f$, where the matrix 
$B$ has been defined by (\ref{a4}). Since ${\tilde W}B{\tilde W}^\ast$ and 
$A_f$ are both self-adjoint operators, the latter inclusion implies 
that ${\tilde W}B{\tilde W}^\ast = A_f$. On the other hand, by Remark 2 
the diagonal matrix $W$ with diagonal entries $w_1$ and 
$w_2$ also satisfies the relation $WBW^\ast =A_f$. Consequently, we have 
$VBV^\ast=B$, where $V=\tilde{W} W^{-1}$. 
Set $v_j:= \tilde{w}_j\overline{w_j}$ for $j=1,2$. Then we obtain 
$VB^2V^\ast=B^2$ which means that 
$v_j e^{4\alpha P} \overline{v_j}=e^{4\alpha P}, j=1,2$. From this 
we conclude that $v_j e^{\ii t 4\alpha P} \overline{v_j}=e^{\ii t 4\alpha P},~
t\in\bbbr$. Thus the function $v_j(x)$ commutes with the translation group 
on the real line. Therefore, $v_j(x)$ is constant a.e. on $\bbbr$, say 
$v_j(x)=\gamma_j$. Then, $\tilde{w}_j=\gamma_j w_j$ for $j=1,2$. 
Inserting this into relation (\ref{b2}), applied to 
$\tilde{w}_j$ and to $w_j$, and remembering that $f\ne 0$ we obtain 
$\gamma_1= \gamma_2=:c$. This completes the proof of Theorem 1. 

Next let us turn to the proof of Theorem 2. Since 
$w_1\in\Hh(2\alpha,0;\varepsilon)$ and 
$w_2\in\Hh(2\alpha,-2\alpha;\varepsilon)$, the 
functions $u_f(z)$ and $g_f(z)$ are indeed holomorphic on the 
strips $\J(2\alpha, 0)$ and $\J(\alpha,-\alpha)$, respectively. The 
relation $L_f^\ast=R_f$ has been proved by Lemma 3(iii).

We verify that $|L_f|=L_{g_f}$. As already noted above, 
we have $L_f=f(x{-}\alpha\ii) e^{2\alpha P}$ by Lemma 3(i). 
The first relation of (\ref{a11}) implies that ${\overline w_1} L_f = 
e^{2 \alpha P} {\overline w_2}$. Using the 
first equality of (\ref{a12}) and the preceding facts we conclude that
\begin{equation}\label{a18}
|L_f| =w_2 e^{2\alpha P}\overline{w_2} 
=w_2\overline{w_1} L_f = w_2\overline{w_1} f(x{-}\alpha\ii)e^{2\alpha P}.
\end{equation}
On the other hand, since 
$\overline{w_1}(x) f(x{-}\alpha\ii)w_2(x{-}2\alpha\ii)=1$ a.e. on $\bbbr$ 
by (\ref{b2}) and $\overline{w_2} (x{-}2\alpha\ii)=w_2(x{-}2\alpha\ii)^{-1}$, 
we obtain
$$
w_2(x) \overline{w_1}(x) f(x{-}\alpha\ii)=w_2(x)\overline{w_2} (x{-}2\alpha\ii)= g_f(x{-}\alpha\ii).
$$
Inserting this into (\ref{a18}) it follows that $|L_f|=L_{g_f}$. 
By (\ref{a18}) and the definition of the function $u_f$ we have 
$L_f= {\overline w_1} w_2 |L_f| = u_f |L_f|$. Therefore, 
the relation $L_f = u_f L_{g_f}$ must be the polar decomposition of the 
closed operator $L_f$. This finishes the proof of Theorem 2.

\bn
\section*{4\ The Example ${\bf f(z)=z}$}

\sn
In this brief section we treat the simplest case $f(z)=z$ and express the 
corresponding functions $w_1$ and $w_2$ in terms of the Weierstra\ss\ Delta function (see,  for instance, [4])
\[
\Delta (z)=z e^{cz}\overset{\infty}{\underset{n=1}{\Pi}} \left( 
1+\frac{z}{n}\right) 
e^{-\frac{z}{n}}~.
\]
Here $c={\lim\limits_{n\rightarrow\infty}} (1+\frac{1}{2} 
+{\dots}+\frac{1}{n} - \log(n{+}1))$ denotes the Euler-Mascheroni constant. 
It is well-known  that $\Delta(z)$ is an entire function which satisfies the 
functional equation
\begin{equation}\label{a27}
\Delta(z)=z\Delta(z+1)~, z\in\bbbc~.
\end{equation}
Let us abbreviate
\begin{equation}\label{a28}
\beta:=\frac{-\ii}{4\alpha}~,~~ \gamma=-\frac{1}{2\alpha} \log 4\alpha
\end{equation}
and define two meromorphic functions $w_1$ and $w_2$ by
\begin{equation}\label{a29}
w_1=e^{\ii\gamma z}\frac{\Delta(\beta z{+}\frac{1}{4})}
{\Delta({-}\beta z{+}\frac{1}{4})},~~w_2 (z)= \ii e^{\ii\gamma z}
\frac{\Delta (\beta z{+}\frac{3}{4})}{\Delta ({-}\beta z{+}\frac{3}{4})}~.
\end{equation}
From the definition of the Delta function it is clear that 
$|w_1(x)|=|w_2(x)|=1$ for real $x$. Further, one easily verifies 
that $w_1\in\Hh (2\alpha, 0;\varepsilon)$ and $w_2\in\Hh(2\alpha, 
-2\alpha;\varepsilon)$ for large $\varepsilon$. Using 
formulas (\ref{a28}) and (\ref{a29}) we conclude that
\begin{align*}
\frac{w_1(z)}{w_2(z{-}2\alpha \ii)}& = \frac{1}{\ii e^{2\alpha\gamma}}~
\frac{\Delta(\beta z{+}\frac{1}{4})}{\Delta({-}\beta z{+}\frac{1}{4})}~
\frac{\Delta({-}\beta(z{-}2\alpha \ii)
{+}\frac{3}{4})}{\Delta(\beta(z{-}2\alpha \ii){+}\frac{3}{4})}\\
&=-4\alpha \ii \frac{\Delta(\beta z {+}\frac{1}{4})}{\Delta(\beta 
z{+}\frac{5}{4})}~=-4\alpha \ii \left(\beta z{+}\frac{1}{4}\right) 
=z-\alpha \ii,
\end{align*}

\begin{align*}
\frac{w_2(z)}{w_1(z{-}2\alpha \ii)} &=\frac{\ii}{e^{2\alpha\gamma}}~
\frac{\Delta(\beta z{+}
\frac{3}{4})}{\Delta({-}\beta z{+}\frac{3}{4})}~
\frac{\Delta({-}\beta (z{-}2\alpha \ii){+}\frac{1}{4})}
{\Delta(\beta(z{-}2\alpha \ii){+}\frac{1}{4})}\\
&=4\alpha \ii\frac{\Delta({-}\beta z{-}\frac{1}{4})}{\Delta 
({-}\beta z {+}\frac{3}{4})}= 4\alpha \ii 
\left({-}\beta z{-}\frac{1}{4}\right) =z{-}\alpha \ii.
\end{align*}
Therefore, by the uniqueness assertion of Theorem 1, $w_1(z)$ and 
$w_2(z)$ are the functions $w_1$ and $w_2$ for the holomorphic 
function $f(z)=z$. In particular we see that $w_1(z)$ has zeros at 
$(4n{+}1)\alpha \ii,~ n\in\bbbn_0$, and poles at 
$-(4n{+}1)\alpha \ii,~ n\in\bbbn_0$, while $w_2(z)$ has zeros at 
$(4n{+}3)\alpha \ii,~ n\in\bbbn_0$, and poles at 
$-(4 n{+}3)\alpha \ii,~ n\in\bbbn_0$. All these zeros and poles are simple and 
there are no other zeros and poles of $w_1$ and $w_2$. 
Inserting $w_1$ and $w_2$ into (\ref{b6}) we obtain explicit 
expressions for the components $u_f$ and $g_f$ of the polar decomposition 
of $f(z)=z$. We omit the details. 

\bn
\section*{5\ An Operator Representation of a $q$-Deformed \break 
Heisenberg Algebra}

\sn
As a by product of the preceding considerations we obtain an interesting 
operator representation of the $q$-deformed Heisenberg algebra introduced 
in [1]. First let us recall the definition of this algebra.

For a positive real number $q\ne 1$, let $\A(q)$ be the complex unital 
algebra with generators ${\bf p},{\bf x},{\bf u},{\bf u}^{-1}$ and 
defining relations 
\begin{eqnarray}\label{a30}
{\bf u}{\bf p}{\bf u}^{-1}=q{\bf p},~~ 
{\bf u}{\bf x}{\bf u}^{-1}=q^{-1}{\bf x},~~ 
{\bf u}{\bf u}^{-1}={\bf u}^{-1} {\bf u}=1,\\
\label{a31}
{\bf p}{\bf x}-q{\bf x}{\bf p}=\ii q^{1/2}(q-q^{-1}){\bf u},~~ 
{\bf x}{\bf p}-q{\bf p}{\bf x}=-\ii q^{1/2}(q-q^{-1}){\bf u}^{-1}.
\end{eqnarray}
If we replace (\ref{a31}) by the relations
\begin{equation}\label{a31a}
{\bf p}{\bf x}=\ii q^{1/2} {\bf u}^{-1}-\ii q^{-1/2} {\bf u},~~ 
{\bf x}{\bf p}=\ii q^{-1/2} {\bf u}^{-1} -\ii q^{1/2} {\bf u},
\end{equation}
then we obtain an equivalent set of defining relations 
(\ref{a30}) and (\ref{a31a}). The 
algebra $\A(q)$ becomes a 
$\ast$-algebra with respect to the involution determined on the generators by
$$
{\bf p}={\bf p}^\ast, {\bf x}={\bf x}^\ast, {\bf u}^\ast={\bf u}^{-1}.
$$
Let us write $q=e^{2\gamma}$ with $\gamma$ real and fix two real 
numbers $\alpha$ and $\beta$ such that $\alpha\beta=\gamma$. In order 
to be in accordance with the preceding sections, one may assume that 
$\alpha > 0$ as well. 

Representations of the $\ast$-algebra $\A(q)$ by Hilbert space operators 
have been constructed and studied in [1] and [7]. We now give an 
operator representation of $\A(q)$ in terms of the "usual" position and 
momentum operators $P=-i{\frac{d}{dx}}$ and $Q=x$ on $L^2(\bbbr)$ 
such that ${\bf p}$ and ${\bf u}$ are 
both represented by self-adjoint operators. 
On the Hilbert space $\Hh=L^2(\bbbr)\oplus L^2(\bbbr)$, we define three 
operators $\rho ({\bf u}), \rho({\bf x})$ by the operator matrices
$$
\rho ({\bf u})=\left( \begin{matrix} e^{\ii \beta x} &0\\ 0 &e^{\ii \beta x}
\end{matrix}\right),~ \rho ({\bf p})=\left(\begin{matrix} 0 &e^{-2\alpha P}\\ 
e^{-2\alpha P} &0
\end{matrix}\right),
$$
$$
\rho ({\bf x})=\left( \begin{matrix}0 &2\sin\beta (x{-}\alpha \ii) 
e^{2\alpha P}\\ 
e^{2\alpha P} 2\sin\beta(x{+}\alpha \ii) &0\end{matrix}\right).
$$
Obviously, $\rho ({\bf u})$ is a unitary and $\rho({\bf p})$ is an unbounded 
self-adjoint operator. Since $\inf~\{|\sin\beta(x{-}\alpha \ii)|;~ 
x\in\bbbr\}>0$, the operator $2\sin\beta(x{-}\alpha i)e^{2\alpha P}$ 
is closed by Lemma 3(i) and so $2\sin\beta(x{-}\alpha \ii)e^{2\alpha P}=
L_{2\sin\beta x}$. Further, since the multiplication operator by the 
function $2\sin\beta (x{-}\alpha \ii)$ is bounded, we have 
$(L_{2\sin\beta x})^\ast= e^{2\alpha P} 2\sin\beta(x{+}
\alpha \ii)$. (Note that because of the boundedness of the function 
$2\sin\beta (x{-}\alpha \ii)$ we do not need the full strength of 
Lemma 3(iii) here.) 
Therefore, $\rho({\bf x})$ coincides with the operator 
$A_{2\sin\beta x}$ defined by (\ref{af}). In particular we conclude in this 
manner that $\rho({\bf x})$ is a self-adjoint operator. 

Next we check that the 
operators $\rho({\bf u}),\rho({\bf p})$ and $\rho({\bf x})$ satisfy the 
relations (\ref{a30}) 
and (\ref{a31a}). Let $\eta\in\D(e^{2\alpha P})$. Then, obviously 
$e^{2\alpha P}\eta\in\D(e^{-2\alpha P})$. From the characterization 
of the domain $\D(e^{-2\alpha P})$ given by Lemma 1 (more precisely, 
from the coresponding assertion for $\alpha<0$) it follows that 
$2\sin\beta (x{-}\alpha \ii) e^{2\alpha P}\eta$ is also 
in $\D(e^{-2\alpha P})$. By (\ref{c10}), we have
 $e^{-2\alpha P}2\sin\beta(x{-}\alpha \ii)
e^{2\alpha P}=2\sin\beta(x{+}\alpha \ii)\eta$. Thus we shown that 
$$
e^{-2\alpha P}2\sin\beta(x{-}\alpha \ii)e^{2\alpha P}\subseteq 
2\sin\beta(x{+}\alpha \ii).
$$
Using this equation and the fact that $q^{\pm 1/2}=q^{\pm\alpha\beta}$ 
we easily get
$$
\rho({\bf p})\rho(\bf x)\subseteq 
\ii q^{1/2} \rho({\bf u})^{-1}{-}\ii q^{-1/2}
 \rho({\bf u}).
$$
In a similar manner we conclude that
\begin{equation*}
\begin{split}
\rho({\bf x})\rho({\bf p})\subseteq \ii q^{-1/2}\rho({\bf u})^{-1}{-}
\ii q^{1/2}\rho({\bf u}),\\
\rho({\bf u})\rho({\bf p})\rho({\bf u})^{-1}=q\rho({\bf p}),~~
\rho({\bf u})\rho({\bf x})\rho{{\bf u}}^{-1}=q^{-1}\rho({\bf x})
\end{split}
\end{equation*}
That is, the operators $\rho({\bf u}), \rho({\bf p})$ and $\rho({\bf x})$ fulfill the 
defining relations (\ref{a30}) and (\ref{a31a}) of the algebra $\A(q)$. 
There is also an invariant dense core for these operators; for instance, 
one may take the domain $\D_0 \oplus\D_0$, where $\D_0$ is 
defined by (\ref{c9}). Thus, we get indeed a 
$\ast$-representation [8] of the $\ast$-algebra $\A(q)$ on the invariant 
dense domain $\D_0 \oplus\D_0$.

Now let $w_1$ and $w_2$ be the functions 
from Theorem 1 in the case $f(z)=2\sin\beta z$. Let $W$ and $V$ denote 
unitary diagonal operator on the Hilbert space 
$\Hh=L^2(\bbbr)\oplus L^2(\bbbr)$ defined by 
$W(\eta_1,\eta_2)=(w_1\eta_1, w_2\eta_2)$ and 
$V(\eta_1,\eta_2)=(V_-\eta_1,V_-\eta_2)$ for 
$\eta_1,\eta_2\in L^2(\bbbr)$, where $(V_-\eta_j)(x):=\eta_j(-x),j=1,2$. 
From formula (\ref{a13}) (or equivalently (\ref{a4})) 
in the proof of Theorem 1 we then obtain
\begin{gather*}
(WV)^\ast \rho ({\bf x}) WV=V^\ast W^\ast A_{2\sin\beta x} WV=
V^\ast W^\ast \left( \begin{matrix}0&L_{2\sin\beta x}\\ 
(L_{2\sin\beta x})^\ast &0\end{matrix}\right) WV\\
V^\ast\left(\begin{matrix} 0&e^{2\alpha P}\\ e^{2\alpha P} 
&0\end{matrix}\right)~V=\left(\begin{matrix} 0 &e^{-2 \alpha P}\\ 
e^{-2 \alpha P} &0\end{matrix}\right) =\rho ({\bf p}).
\end{gather*}
That is, the self-adjoint operators $\rho({\bf x})$ and $\rho({\bf p})$ are 
unitarily equivalent via the unitary operator $WV$. Thus, if we 
think of $\rho({\bf x})$ and $\rho({\bf p})$ as $q$-analogues of the 
position and momentum operators, then the unitary $WV$ takes the role of the 
Fourier transform in these considerations. Being the main ingredients of 
the unitary operator $WV$ the two holomorphic functions $w_1$ and $w_2$ 
from Theorem 1 are crucial in this context.

\bn

\vspace*{12pt}
\centerline{Konrad Schm\"udgen}
\centerline{Fakult\"at f\"ur Mathematik und Informatik}
\centerline{Universit\"at Leipzig, Augustusplatz 10, 04109 Leipzig, Germany}
\centerline{E-mail: schmuedg@mathematik.uni-leipzig.de}

\end{document}